\documentclass[11pt]{amsart}  
\usepackage{amsmath,amssymb,amsfonts,graphicx}
\author[Christine Bachoc, Yael Ben-Haim, and Simon Litsyn]{Christine Bachoc, Yael Ben-Haim,\\ and Simon Litsyn, {\em Senior Member, IEEE} }
\thanks{C. Bachoc is with Laboratoire A2X,
Universit\'e Bordeaux I, 351, cours de la Li\-b\'e\-ration, 33405 Talence
France. (email: bachoc@math.u-bordeaux1.fr)
}
\thanks{Y. Ben-Haim and S. Litsyn are with the School of Electrical Engineering, Tel-Aviv University, Tel-Aviv 69978, Israel (email: yael@eng.tau.ac.il; litsyn@eng.tau.ac.il).
}

\title[Bounds for codes in products of spaces]{Bounds for Codes in Products of Spaces, Grassmann and Stiefel Manifolds}

\date{\today}

\newtheorem{defi}{Definition}[section]
\newtheorem{definition}[defi]{Definition}
\newtheorem{proposition}[defi]{Proposition}
\newtheorem{theorem}[defi]{Theorem}
\newtheorem{remark}[defi]{Remark}

\newtheorem{lemma}[defi]{Lemma}

\newcommand{\N}{{\mathbb{N}}} 

\newcommand{\R}{{\mathbb{R}}} 
\newcommand{\C}{{\mathbb{C}}}
\newcommand{\HH}{{\mathbb{H}}}
\newcommand{\OO}{{\mathbb{O}}}

\newcommand{\Gmn}{{\mathcal{G}_{m,n}}}
\newcommand{\Vmn}{{\mathcal{V}_{m,n}}}

\newcommand{\uk}{\underline{k}}
\newcommand{\ul}{\underline{l}}
\newcommand{\Sn}{S^{n-1}}

\newcommand{\Snm}{\big(S^{n-1}\big)^m}
\newcommand{\Scnm}{\big(S^{cn-1}\big)^m}
\newcommand{\Pn}{{\mathbb P}^{n-1}}
\newcommand{\PO}{{\mathbb P}^{2}(\OO)}
\newcommand{\Pnm}{{\big(\mathbb P}^{n-1}\big)^m}

\newcommand{\Pk}{P_{\uk}}

\newcommand{\Pl}{P_{\ul}}

\newcommand{\Harm}{\operatorname{Harm}}

\newcommand{\Id}{\operatorname{Id}}

\newcommand{\trace}{\operatorname{trace}}

\newcommand{\te}{\theta}
\newcommand{\vp}{\varphi}
\newcommand{\al}{\alpha}

\begin{document}

\begin{abstract} 
Upper bounds are derived for codes in Stiefel and
Grassmann manifolds with given minimum chordal distance. They stem 
from upper bounds for codes in the product of unit spheres and
projective spaces. The new bounds are asymptotically better than 
the previously known ones.
\end{abstract}

\maketitle

\section{Introduction}  \label{sec intro}

Use of multiple transmit and receive antennas essentially increases the spectral efficiency
of wireless systems (see \cite{ZT02} and references therein). Analysis of Rayleigh flat-fading multiple-input multiple-output (MIMO) 
scenarios with $m$ transmit antennas and $n$ transmitted symbols, reveals that relevant
coding schemes can be designed as collections of elements (points) in
the complex Grassmann manifold - the set of $m$-dimensional linear subspaces in $\C^n$, 
if the channel is unknown to the receiver, and in the complex Stiefel manifold - the set of $m$ orthonormal 
vectors in $\C^n$, if the channel is known to the receiver. An appropriately defined distance measure 
between the points characterizes diversity of the designed scheme. 
Following standard for coding theory considerations,
we study the relation between the number of points (the size of a code) and the minimum distance between 
distinct code points. Our aim in this paper is to obtain new upper bounds for the size of
codes in Grassmann and Stiefel manifolds.  

The most powerful technique for this kind of problems is the linear programming method (called also the polynomial method), 
initiated by Delsarte \cite{Del73}. The method is very
well understood in the case of $2$-point homogeneous spaces (defined in the next section), where
very explicit bounds, and also good asymptotic bounds on the rate
of codes have been derived. Examples are the Hamming and Johnson schemes, treated in \cite{MRRW}, 
the unit
sphere
of $\R^n$ \cite{KL78}, and the projective spaces \cite{KL78}, \cite{Lev98}.

When the underlying space is homogeneous and  symmetric but 
not $2$-point homogeneous, the situation
is much more complicated, although the principles of the linear
programming method remain valid. The difficulties come from the fact
that the zonal functions defined for these spaces are not functions of one variable, but
afford several variables. The Grassmann spaces considered in this
paper fall into this category. An attempt to overcome this problem was
carried out in \cite{Bac06}. An asymptotic bound for the rate of
Grassmannian codes was
obtained, involving the asymptotics of the largest eigenvalue of some
symmetric endomorphism. This bound however is not optimal since it was
improved for $m>1$ by some volume-type arguments for a large range of values of the
minimal distance \cite{BN06}. 

There is one trivial case of symmetric spaces of rank $m>1$ for which
the classical treatment of the linear programming method is easily extended: it is  the direct product of $2$-point homogeneous spaces,
such as the direct product of $m$ copies of the unit sphere. An example of a similarly easy case is provided by the non-binary Johnson space \cite{Aal90,TAG,BL06}, that is the product of the Hamming and the binary Johnson scheme.

The approach developed this paper is to relate Grassmann and Stiefel
spaces  and their associated codes to various products of 2-point
homogeneous spaces,
and hence to derive upper bounds for these codes in a rather easy way. 
The asymptotic versions
of the new bounds provide the best currently known asymptotic bounds.

The paper is organized as follows. Definitions and known results are
given in Section \ref{sec def}. 
Section \ref{sec relations} describes various relations between
the spaces and their codes. The  simplest one connects Grassmann
and Stiefel
spaces to the unit sphere of an asymptotically equal dimension;
therefrom, for example, we obtain a bound for the asymptotic rate of Grassmannian
codes that already improves upon the previous ones (Theorem \ref{bound
grass}).
Section \ref{sec products} develops the Delsarte polynomial
method for the products of spaces under consideration,
including the classical method that involves the Christoffel-Darboux
formula, and derives upper bounds for the size of the associated
codes. A bound for the asymptotic rate of these codes is obtained. 
Section \ref{sec grass stief} discusses the consequences for the
Grassmannian and Stiefel codes. In particular, we show that the bound
obtained on the asymptotic rate of Grassmannian codes from the product
of projective spaces is sometimes better than the one obtained in Section \ref{sec relations}. We conclude in Section \ref{sec con}.

\section{Definitions and known results}     \label{sec def}

We shall use the following notations and definitions. We say that $f(n)\lesssim g(n),~f(n)\simeq g(n),~f(n)\gtrsim g(n)$ if
$\lim_{n\to\infty}\frac{f(n)}{g(n)}\leq1,~\lim_{n\to\infty}\frac{f(n)}{g(n)}=1,~\lim_{n\to\infty}\frac{f(n)}{g(n)}\geq1$,
respectively. A code in a metric space $(X,d)$ is a finite set contained in the space, and a
codeword is an element of the code. The size of a code $C$ is its
cardinality, and its rate is $R(C):=\frac{1}{n}\ln|C|$ where $\ln$
denotes the natural logarithm. The meaning of
$n$ will be defined separately for each space. Indeed, a
more consistent and general definition of the rate of a code in a
manifold $X$ would be $\frac{1}{\dim(X)}{\ln|C|}$. The minimum distance of a code is the minimum distance (induced by the relevant metric) between a pair of distinct codewords.
A metric space $(X,d)$ is called 2-point homogeneous, if $X$ affords
the transitive action of a group $G$, such that the orbits of the
action of $G$ on $X\times X$ are characterized by the distance $d$. In
other words, for all $(x,y)\in X$ and $g\in G$, $d(g(x),g(y))=d(x,y)$,
and moreover, for all pairs $(x,y)$, $(x',y')\in G$, there exists
$g\in G$ such that $g(x)=x'$ and $g(y)=y'$ if and only if
$d(x,y)=d(x',y')$. It
is a well-known fact that the compact Riemannian manifolds that are two-point
homogeneous are exa\-ctly: the unit sphere $\Sn$, the projective spaces
$\Pn(K)$ where $K=\R,\C,\HH$ and the projective plane over the
octonions $\PO$ (see \cite{Hel78}, and \cite{CS03} for more about the
octonions and $\PO$).

\subsection{The real compact two-point homogeneous spaces}
The unit \\ sphe\-re of the Euclidean space $\R^n$ is denoted
$\Sn$, namely,
\begin{equation}
\Sn := \left\{ (x_1,\ldots,x_n)\in\R^n~|~\sum_{i=1}^n x_i^2=1\right\}.
\end{equation}
The standard scalar product in $\R^n$, given by 
$(u\cdot v)=\sum_{i=1}^n u_iv_i$,
defines the Euclidean distance between two points of $\Sn$:
\begin{equation}
\Vert u-v\Vert =\sqrt{\sum_{i=1}^n (u_i-v_i)^2}
=\sqrt{2}\sqrt{1-(u\cdot v)}.
\end{equation}
The angular distance between $u$ and $v$ is defined by the angle
$\theta(u,v)\in [0,\pi]$, also denoted $\theta$ when the context is clear.
We have of course
\begin{equation}
\Vert u-v\Vert =\sqrt{2}\sqrt{1-\cos\theta}.
\end{equation}

The best known asymptotic bounds on the rate of spherical codes as a
function of the minimum distance are given 
in the following theorems. Theorem \ref{theo gv sphere} is the Chabauty-Shannon-Wyner, which is the analog of the Gilbert-Varshamov
bound for codes over finite fields, Theorem \ref{theo lp sphere} is the linear programming bound, and Theorem \ref{theo yag sphere} is an Elias-Bassalygo type improvement due to Yaglom. See \cite{EZ01,CS88} for details on these bounds.

\begin{theorem}[\cite{Sha59},\cite{Cha53},\cite{Wyn65}] \label{theo gv sphere}
There exists a sequence of codes $\{C_n\}_{n=1}^\infty$ such that $C_n$ is of length $n$, minimum angular distance $\theta\leq\pi/2$, and rate $R(C_n):=\frac{1}{n}\ln|C_n|$ which satisfies
\begin{equation}
R(C_n) \gtrsim -\ln\sin\theta.
\end{equation}

\end{theorem}

\begin{theorem}[\cite{KL78}] \label{theo lp sphere}
Let $C$ be a spherical code with minimum angular distance $\theta\leq\pi/2$. Then, when $n\to\infty$,
\begin{equation}    \label{eq kl}
R(C) \lesssim R_{LP}(\te) := \frac{1+\sin\te}{2\sin\te}H\left(\frac{1-\sin\te}{1+\sin\te}\right),
\end{equation}
where $H$ is the entropy function, $H(x):=-x\ln x-(1-x)\ln(1-x)$.

\end{theorem}

\begin{theorem}[\cite{KL78}]    \label{theo yag sphere}
Let $C$ be a spherical code with minimum angular distance $\theta\leq\pi/2$. Then, for any $\te/2\leq\vp\leq\pi/2$, when $n\to\infty$,
\begin{equation}    \label{eq yag sphere}
R(C) \lesssim \tilde{R}(\al) -\ln\sin\vp,
\end{equation}
where $\sin(\al/2)=\frac{\sin(\te/2)}{\sin\vp}$, and $\tilde{R}(\al)$ is an upper bound on the asymptotic rate of spherical codes with minimum angular distance $\al$.

\end{theorem}
For $\tilde{R}=R_{LP}$, the optimal choice of $\varphi$ is discussed in \cite{KL78}; it
corresponds to $\alpha\sim 63^{\circ}$ and gives a better bound than 
$R_{LP}(\theta)$ when $\theta$ is smaller than $\alpha$:
\begin{equation}\label{Y}
R(C)\lesssim R_Y(\theta):=-\ln\sqrt{1-\cos\theta}-0.0686.
\end{equation}

We denote by $R_S(\theta)$ the function that provides the best known bound for the asymptotic
rate of spherical codes:
\begin{equation}\label{RS}
R_S(\theta)=\left\{
\begin{array}{ll}
R_Y(\theta) &\text{ if } \theta<\alpha\\
R_{LP}(\theta) &\text{ if } \alpha\leq \theta \leq \pi/2
\end{array}
\right.
\end{equation}

\begin{remark}
The asymptotic rate of spherical codes with minimum angular distance at least $\pi/2$ is known to be equal to zero. This is a consequence of the Rankin bound (\cite{Ran55}, see also \cite{EZ01} or \cite{CS88}). 
\end{remark}

The other real compact manifolds which are two-point homogeneous can
be treated in a similar way. As was recalled before, these are the projective
spaces $\Pn(K)$ where $K=\R,\C,\HH$ (the field of real quaternions)
and $n\geq 3$, and the projective plane over the
octonions $\PO$. 
In order to treat the fields of coefficients 
in a uniform way, we extend the definition of $(x\cdot y)$ so that,
for all $x,y\in K^n$,
$(x\cdot y)=\sum_{i=1}^n x_i \overline{y_i}$, where the conjugation $x\to
\overline{x}$ is the standard one over $K=\C,\HH,\OO$ and is the
identity over $\R$. Also we conventionally
assume that $n=3$ when $K=\OO$. The group $G$ under which these
spaces are two-point homogeneous is respectively the orthogonal group
$O(\R^n)$, the unitary groups $U(K^n)$ with $K=\C,\HH$, and the Lie
group $F_4$ (see \cite
{CS03} for this last case).

The angular distance between $p$ and $q$ in $\Pn(K)$ is defined by the
angle $\theta=\theta(p,q)\in [0,\pi/2]$ such that
$\cos\theta(p,q)=|(e\cdot f)|$ 
where $e$, $f$ are arbitrary chosen unit vectors of the lines $p$,
$q$. It is shown in \cite{KL78} and \cite{Lev98} that the linear programming method
applies to these spaces. The derived asymptotic bound for the rate can 
also be obtained from the bounds for spherical codes, because to a
code $C$ in $\Pn(K)$ one can associate a code in $S^{cn-1}$ with the same
size and a minimum angular distance at least equal to the one of $C$, 
selecting a unit vector in each element of
$C$. One obtains:

\begin{theorem}[\cite{KL78}] \label{theo lp proj}
Let $C$ be a code in $\Pn(K)$  with minimum angular distance
$\theta\leq\pi/2$. 
Let $c:=1,2,4$ respectively when $K=\R,\C,\HH$ (so that $c=[K:\R]$). Then, when $n\to\infty$,
\begin{equation}    
R(C):=\frac{1}{n}\ln |C| \lesssim cR_{S}(\te) 
\end{equation}
\end{theorem}

\subsection{The Grassmann space}     \label{sec def grass}

Let $K$ be the real or the complex field. The Grassmann space
$\Gmn(K)$ is the set of all subspaces of dimension $m$ in $K^n$. It is
a homogeneous space under the action of either the orthogonal group
$O(\R^n)$ or the unitary group $U(\C^n)$. We will denote $\Gmn$ when
$K$ is arbitrary. It is worth noticing that when $m=1$ we recover the
projective space. Several metrics have been defined in $\Gmn$, see
\cite{BN02,EAS}. In this paper we consider the {\em chordal distance},
which was introduced in \cite{CHS} and studied in \cite{BN02,BN06,Bac06,CHS,EAS}. The following two definitions for the chordal distance $d_c(p,q)$ are equivalent.

\begin{definition}  \label{def chord1}
Given the planes $p,q\in\Gmn$, apply the following procedure. Initialize the sets of unit vectors $A=\emptyset$ and $B=\emptyset$. In the $i$th step, choose the vectors $a_i,b_i$ such that:
\begin{enumerate}
\item $a_i$ is contained in $p$ and $b_i$ is contained in $q$.
\item $a_i$ is orthogonal to all the vectors in $A$ and $b_i$ is orthogonal to all the vectors in $B$.
\item Among all the vectors satisfying the conditions in (i) and (ii),
  the angle between $a_i$ and $b_i$ is minimal (i.e., their inner /
  Hermitian product module is maximal).
\end{enumerate}
Set $\te_i$ to be the angle between $a_i$ and $b_i$, insert $a_i$ to
$A$ and $b_i$ to $B$, and proceed to the next step, until $m$ angles
$0\leq\te_1\leq\te_2\leq\ldots\leq\te_m\leq \pi/2$, called the
principal angles between $p$ and $q$, have been defined. Then the chordal distance is
\begin{equation*}
d_c(p,q):=\sqrt{\sum_{i=1}^m \sin^2\theta_i}=\sqrt{m-\sum_{i=1}^m\cos^2\te_i}.
\end{equation*}

\end{definition}


\begin{lemma}[\cite{CHS}]
For a plane $p\in\Gmn$, let $A_p$ be a $p\times n$ matrix whose rows
form an orthonormal basis of $p$, 
and let $\pi_p:=A_p^*A_p$ be the matrix of the orthogonal projection
on $p$ ($A_p^*$ denotes the Hermitian conjugate of $A_p$). Then, the
projection matrix $\pi_p$ does not depend on the choice of $A_p$, and,
given two planes $p,q\in\Gmn$, the chordal distance is
\begin{equation}    \label{eq minchord2}
d_c(p,q):=\sqrt{m-\trace(\pi_p\circ\pi_q)}.
\end{equation}

\end{lemma}

We review some known bounds on the size of codes in Grassmann
spaces. We recall that $c=1$ if $K=\R$ and $c=2$ if $K=\C$. The first work in the area is due to Conway et al.
\cite{CHS}. They found an isometric embedding from $\Gmn(\R)$ to the Euclidean sphere with radius $\sqrt{\frac{m(n-m)}{n}}$ in $\R^{\frac{1}{2}(n-1)(n+2)}$ (a similar embedding exists also for $\Gmn(\C)$, see \cite{BN02}). This enabled them to use the Rankin bound on spherical
codes \cite{Ran55} to derive the following non-asymptotic bound.
\begin{theorem}
Let $C$ be a code in $\Gmn(\R)$ with minimum chordal distance $d$. Then
\begin{equation}
d \leq \left\{ \begin{array}{cc} \frac{m(n-m)}{n}\frac{|C|}{|C|-1} &    \\
\frac{m(n-m)}{n} & \text{if } |C|>\frac{1}{2}n(n+1)
\end{array} \right.
\end{equation}

\end{theorem}

Later, an asymptotic expression to the volume of a ball in Grassmann spaces was derived by Barg and Nogin \cite{BN02},
yielding the analogue of Gilbert-Varshamov and the Hamming asymptotic  bounds for
the rate $R(C):=\frac{1}{n}\ln|C|$ of Grassmannian codes.
\begin{theorem}[\cite{BN02}]
For any pair of constants $d,m,$ such that $d\leq\sqrt{m}$, there exists an infinite sequence of codes $\{C_n\}$ in $\Gmn$ with minimum chordal distance $d$ and rate
\begin{equation}
R\gtrsim -cm\ln\frac{d}{\sqrt{m}}.
\end{equation}

\end{theorem}
\begin{theorem}[\cite{BN02}]
Let $C$ be in $\Gmn$ with minimum chordal distance $d$. Then, when $n\to\infty$,
\begin{equation}    \label{eq bn1}
R(C) \lesssim -cm \ln \left( \sqrt{ 1-\sqrt{ 1- d^2/2m} } \right) .  
\end{equation}

\end{theorem}

A linear programming bound was derived by Bachoc \cite{Bac06}.
\begin{theorem}[\cite{Bac06}]   \label{theo gras lp}
Let $C$ be a code in $\Gmn(\R)$ with minimum chordal distance $d$. Then, when $n\to\infty$,
\begin{equation}    \label{eq gras lp}
R(C) \lesssim m[(1+\rho)\ln(1+\rho)-\rho\ln \rho],
\end{equation}
where
\begin{equation}
\rho = \frac{1}{2}m(\sqrt{m}/d-1).
\end{equation}

\end{theorem}
We note that the derivation of the linear programming bound is not a
straightforward analogy to the derivation of this bound in other
metric spaces, since it involves a family of orthogonal generalized
Jacobi polynomials with several variables, and that the bound (\ref{eq
  gras lp}) in the case $m=1$ coincides with the bound for the real
projective space of Theorem \ref{theo lp proj}. The case of the complex
Grassmann space is not treated in \cite{Bac06} but could be treated
in a similar way.

Finally, a new upper bound was introduced recently by Barg and Nogin \cite{BN06}, using Blichfeldt's density method.
\begin{theorem}     \label{theo blich}
Let $C$ be a code over the Grassmann space $\Gmn(\R)$ with minimum
chordal distance $d$. Then
\begin{equation}    \label{eq blich}
R(C) \lesssim -m\ln\left( \sqrt{ 1-\sqrt{ 1- d^2/m} } \right).
\end{equation}

\end{theorem}
It is immediate to see that this bound improves upon (\ref{eq bn1})
and upon (\ref{eq gras lp}) for $m>1$ and for a large range of  values of $d$. Hence, until this paper, Theorems \ref{theo gras lp} and
\ref{theo blich} provide the best asymptotic known upper bounds on the rate of codes in $\Gmn(\R)$.


\subsection{The Stiefel manifold}   \label{sec def stiefel}

The Stiefel manifold $\Vmn(K)$ is the set of $m$-tuples of orthonormal vectors in $K^n$, or equivalently
\begin{equation*}
\Vmn(K)=\{X\in M^{m\times n}(K) \mid XX^*=\Id_m\},
\end{equation*}
where $\Id_m$ is the $m\times m$ identity matrix. The orthogonal group
$O(\R^n)$ if $K=\R$, respectively the unitary group $U(\C^n)$ if
$K=\C$  acts transitively on
$\Vmn(K),$ and this space can be identified with the
set of classes $O(\R^n)/O(\R^{n-m})$, respectively $U(\C^n)/U(\C^{n-m})$.

The distance considered in coding theory is 
\begin{equation*}
d(X,Y):=\| X-Y \| = \sqrt{\trace((X-Y)(X^*-Y^*))}.
\end{equation*}
In other words, $d(X,Y)$ is the Euclidean distance between $X$ and $Y$, when $X$ and $Y$ are regarded as one-dimensional vectors of length $mn$. We refer the reader to \cite{Hen05} for a treatment of codes in Stiefel manifolds.

\section{More spaces and their interconnections}\label{sec relations}

The simplest of these connections relate Grassmann and Stiefel
spaces to a single unit sphere, and allow to apply directly the known
bounds for spherical codes to the Grassmannian and Stiefel codes. We
start with them, then we introduce the products of spheres and
projective spaces and their relations with Grassmann and Stiefel spaces.

\subsection{$\Gmn$ and $S^{cmn-1}$}

We follow the notations and
definitions of Section \ref{sec def grass}. For all $p,q\in \Gmn(K)$,
we set
\begin{equation*}
\sigma(p,q):=\sum_{i=1}^m \cos^2\theta_i=\trace(\pi_p\circ \pi_q).
\end{equation*}

We define a mapping
\begin{equation*}
\beta : \Gmn(K)\to   S^{cmn-1}
\end{equation*}
in the following way. We select for all $p\in\Gmn(K)$, an
orthonormal basis $(e_1,\dots, e_m)$ 
of $p$ whose elements belong to $K^n$. With the usual identification of $\C$ and $\R\times \R$ through the mapping
$z=x+iy\to (x,y)$, we consider these elements in $\R^{cn}$.  Then
$\beta(p)$ is chosen to be the element of $\R^{cmn}$ obtained by the
concatenation of $e_1,\dots,e_m$, divided by $\sqrt{m}$. Obviously, $\beta(p)\in S^{cmn-1}$.

The new bounds for Grassmann spaces rely on the following lemma.

\begin{lemma}   \label{lem proj grass}
For all $p,q\in \Gmn(K)$,
\begin{equation*}
\cos\te(\beta(p),\beta(q))\leq \sqrt{\frac{\sigma(p,q)}{m}}.
\end{equation*}
\end{lemma}

\begin{proof}  Let $\beta(p)=e$, obtained from an orthonormal basis
$(e_1,e_2,\dots,e_m)$ of $p$
and $\beta(q)=e'$, obtained from an orthonormal basis
$(e'_1,\dots,e'_m)$ of $q$. We compute $\sigma(p,q)=\trace(\pi_p\circ\pi_q)$.
Let $A_p$, $A_q$ denote the $m\times n$ matrices whose rows
are the basis elements $e_i$, $e'_i$ respectively. Then
\begin{equation*}
\sigma(p,q)=\trace(\pi_p\circ \pi_q)=\trace(A_p^*A_pA_q^*A_q)=\trace(A_pA_q^*A_qA_p^*).
\end{equation*}

The entries of the matrix $A_pA_q^*$ are the hermitian products
$(e_i\cdot e'_j)$. So we obtain:
\begin{equation}\label{sigma}
\sigma(p,q)=\sum_{1\leq i,j\leq m}|(e_i\cdot e'_j)|^2
\end{equation}

If $K=\R$, we obtain from Cauchy-Schwartz inequality
\begin{equation}\label{ineq2}
\cos \theta(e,e')=(e\cdot e')=\frac{\sum_{i=1}^m (e_i\cdot e'_i)}{m} \leq
\sqrt{\frac{\sum_{i=1}^m(e_i\cdot e'_i)^2}{m}}\leq  \sqrt{\frac{\sigma(p,q)}{m}}.
\end{equation}

If $K=\C$, let us denote by $\Re(z)$ the real part of a complex number $z$.
In the identification $\C^n=\R^{2n}$ recalled above, the standard
scalar product on $\R^{2n}$ is given by $\Re(h(x,y))$. With the obvious
inequality ${\Re(h(x,y))}^2\leq |h(x,y)|^2$, we obtain the same
inequality $\cos \theta(e,e')\leq \sqrt{\frac{\sigma(p,q)}{m}}$ (where
$e$ and $e'$ are considered in the unit sphere of $R^{2n}$.)

\end{proof}

Let us recall Definition \ref{def chord1} of the chordal distance
$d_c(p,q)$ in Grassmann spaces. 
The definition involves the construction of orthonormal basis
$(a_1,\ldots,a_m)$ and $(b_1,\ldots,b_m)$ for $p$ and $q$
respectively, such that the principle angles $\theta_i$ satisfy $\cos
\theta_i= |(a_i\cdot b_i)|$ and
$\sigma(p,q)=\sum_{i=1}^m\cos^2\te_i$. Lemma \ref{lem proj grass}
shows that if one chooses {\it arbitrary} orthonormal basis
$(e_1,\dots,e_m)$ of $p$, $(e'_1,\dots,e'_m)$ of $q$, and defines an
alternative set of ``principal angles'' 
$\te_1',\dots,\te_m'$ by $\te_i'=\arccos|(e_i\cdot e_i')|$, then
$\sigma(p,q)\geq\sum_{i=1}^m\cos^2\te_i'$. Thus, 
an upper bound on the chordal distance between $p$ and $q$ is obtained.

\smallskip
It follows from Lemma \ref{lem proj grass} that the bounds for spherical codes (of the sphere $S^{cmn-1}$) can
be applied to codes in Grassmann spaces. We obtain for the
asymptotic rate:

\begin{theorem} \label{bound grass} 
Let $C$ be a code in $\Gmn(K)$ with minimal chordal distance $d=\sqrt{m-s}$,
and let $\theta=\arccos(\sqrt{s/m})$. Then, when $n\to +\infty$,
\begin{equation}\label{bound1 gmn}
R(C)\lesssim cmR_S(\te)
\end{equation}
where $R_S$ is defined in \eqref{RS}.

\end{theorem}

\begin{figure}
\begin{center}
\includegraphics*{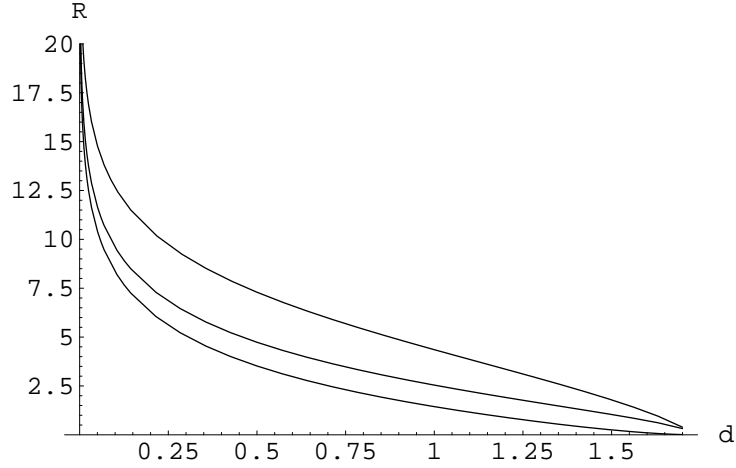}
\caption{Upper bound on the asymptotic rate of real Grassmannian codes with minimum chordal distance $d$, $m=3$. From top to bottom: a linear programming bound \eqref{eq gras lp}, a Blichfeldt-type bound \eqref{eq blich}, and the new bound \eqref{bound1 gmn}}   \label{fig grass}
\end{center}
\end{figure}

Figure \ref{fig grass} compares \eqref{bound1 gmn} with some of the existing bounds, given in Section \ref{sec def grass}. The bound \eqref{eq blich} can also be obtained using
  the mapping $\beta$, joined with  the Rankin-Blichfeldt bound
  \cite{BN06}. It is well known that the Rankin-Blichfeldt bound \cite{Ran55} is
  improved by the bound $R_S(\theta)$ for all values of $\theta$.

\begin{remark}\label{remark}
In \cite{CHS}, the authors introduce 
an isometric embedding of $\Gmn(\R)$ into a unit sphere, but
the dimension of this sphere is much larger than the one of the
Grassmann space ($(n-1)(n+2)/2$ instead of $m(n-m)$). Instead, the
dimension of $S^{cmn-1}$ is asymptotically equivalent to the one of
$\Gmn(K)$, but our embedding is not isometric.

\end{remark}

\subsection{$\Vmn$ and $S^{cmn-1}$}

\begin{lemma}
Let $X,Y\in \Vmn(K)$, $K=\R,\C$.
Let $(e_1,\dots e_m)$ denote the rows of $X$, respectively
$(e'_1,\dots e'_m)$ for the rows of $Y$. Then
\begin{equation*}
d(X,Y)=\sqrt{2}\sqrt{m-\sum_{i=1}^m \Re (e_i\cdot e'_i)}.
\end{equation*}
\end{lemma}

\begin{proof} We calculate
\begin{align*}
\|X-Y\|^2&=\trace((X-Y)(X^*-Y^*))\\
&=\trace(XX^*-XY^*-YX^*+YY^*)\\
&=2m-2\Re(\trace(XY^*))
\end{align*}
since $XX^*=YY^*=\Id_m$ and 
$\overline{\trace(XY^*)}=\trace(\overline{X}\overline{Y^*})=\trace(YX^*)$.
We conclude with
\begin{equation*}
\trace(XY^*)=\sum_{i=1}^m (e_i\cdot e'_i).
\end{equation*}
\end{proof}

Again with the identification of $\C^n$ with $\R^{2n}$, we 
view $\Vmn(\C)$ as a submanifold of ${\mathcal  V}_{m,2n}(\R)$ endowed
with the distance 
$$d(X,Y)=\Vert X-Y\Vert =\sqrt{2}\sqrt{m -\sum_{i=1}^m (e_i\cdot e'_i)}.$$

We see that the obvious mapping:
\begin{align*}
\gamma: \Vmn(K) &\to S^{cmn-1}\\
X &\mapsto \gamma(X)=\frac{1}{\sqrt{m}}(e_1,\dots,e_m)
\end{align*}
is this time, up to a suitable scaling of the distances, an
isometry. Hence the bounds  for spherical codes also apply to $\Vmn(K)$,
probably in a quite efficient way. Still,
one constraint is not encoded in it: the fact that the vectors
$e_i$ are pairwise orthogonal and of norm $1$. We resume these observations in the following theorem.
\begin{theorem}
Any upper bound on the size of codes in $S^{cmn-1}$ with minimum angular distance $\theta=\arccos(s/m)$ is
also an upper bound for codes in the Stiefel space
$\Vmn(K)$ with minimum distance $\sqrt{2}\sqrt{m-s}$. 
In particular, we have for the asymptotic rate, when $n\to +\infty$,
\begin{equation}\label{bound1 stiefel}
R(C)\lesssim cmR_S(\te)
\end{equation}
where $R_S$ is defined in \eqref{RS}.

\end{theorem}

\subsection{$\Gmn$, $\Vmn$ and products of spaces}\label{subsec C}

So far we have established a relation between codes in $\Gmn$ and
$\Vmn$ and codes in $S^{cmn-1}$. It is worth noticing that the
mappings $\beta$ and $\gamma$ defined above factor out by $\Scnm$,
since the elements $e_i$ are unit vectors. Hence  bounds
for codes in $\Scnm$,  will imply bounds for codes in $\Gmn$ and
$\Vmn$. This is the motivation to the  generalization of the
linear programming method to the product of unit spheres, and more
generally to the product of 2-point homogeneous spaces, which is
proposed in the next section. As we shall see, the asymptotic bound for
the rate of codes in $\Scnm$ is not better than for $S^{cmn-1}$, hence
doesn't improve on \eqref{bound1 gmn} and \eqref{bound1 stiefel}. A better result is
obtained for Grassmann spaces with the product of projective spaces.

\smallskip
We now define more precisely the products of spaces and their
associated distances that will be studied in the next section.
We start with the product of $m$ copies of the unit sphere of $\R^n$:
$$\Snm=\{e=(e_1,\dots,e_m)\mid e_i\in \Sn\}.$$
We consider on $\Snm$ the distance given by
\begin{align*}
d(e,e')&=\sqrt{\sum_{i=1}^m \Vert e_i-e'_i \Vert^2}\\
&=\sqrt{2m}\sqrt{1-\frac{\sum_{i=1}^m \cos\theta_i}{m}},
\end{align*}
where $\cos\te_i=(e_i\cdot e_i')$. We attach to a pair $e,e'\in\Snm$ an angle 
$\theta=\theta(e,e')\in [0,\pi]$
such that
\begin{equation}
\cos\theta=\frac{\sum_{i=1}^m \cos\theta_i}{m}
\end{equation}
and we call $\theta$ the {\em angular distance} between $e$ and
$e'$. The angle $\theta$ is also the angle between the vectors
$e/\sqrt{m}$ and $e'/\sqrt{m}$, viewed as elements of $S^{mn-1}$.

\smallskip
We define, for the remaining 2-point homogeneous spaces recalled
above, and without specifying the field $K$,
\begin{align*}
\Pnm=\{p=(p_1,\dots,p_m)\mid p_i\in \Pn\}.
\end{align*}
We attach to a pair $p,p'\in \Pnm$ an angle 
$\theta=\theta(p,p')\in [0,\pi/2]$
such that
\begin{equation*}
\cos^2\theta=\frac{\sum_{i=1}^m \cos^2\theta_i}{m}
\end{equation*}
where $\theta_i=\theta(p_i,q_i)$
and we call $\theta$ the {\em angular distance} between $p$ and $p'$. 
We consider on $\Pnm$ the ``chordal'' distance given by
\begin{equation*}
d(p,p')=\sqrt{m-\sum_{i=1}^m \cos^2\theta_i}=\sqrt{m}\sin\theta.
\end{equation*}

In order to derive bounds for codes in Grassmann spaces
$\Gmn(K)$, we shall make use of the mapping:
\begin{equation*}
\nu:\Gmn(K)\to \big(\Pn(K)\big)^m
\end{equation*}
defined in the following way: for all $p\in \Gmn(K)$, we choose a
$m$-tuple $(p_1,\dots,p_m)$ of pairwise orthogonal lines of $p$. 
We set $\nu(p)=(p_1,\dots,p_m)$. 

Because of the equation \eqref{sigma}, we have similarly:
\begin{equation*}
\cos^2 \theta(\nu(p),\nu(q))\leq \frac{\sigma(p,q)}{m}
\end{equation*}
hence the bounds for codes in $\big(\Pn(K)\big)^m$ apply to codes
in Grassmann spaces.

\section{Bounds for codes in the product of 2-point homogeneous spaces}\label{sec products}

In this section, $X$ denotes one of the spaces $\Sn$,  $\Pn(K)$ where $K=\R,\C,\HH$,
or the projective plane over the
octonions $\PO$. 
We derive bounds for codes in $X^m$ with a given minimum distance,
following Delsarte's linear programming method
as performed in \cite{KL78}. As a reference on orthogonal polynomials, we refer to \cite{Sze39}.

\subsection{Review of the necessary material on the harmonic analysis
  of the spaces $X$}
We recall that, to each of these spaces is associated a family of
orthogonal polynomials of one variable, which are the zonal
polynomials relative to the action of the group $G$ (see \cite{KL78},
\cite{Lev98}, \cite{Hog92}).  For
$X=\Sn$, these polynomials are the Gegenbauer polynomials with
parameter $n/2-1$ and associated orthogonal measure
$(1-x^2)^{(n-3)/2}$ on the interval $[-1,1]$. For $X=\Pn(K)$, these
polynomials are Jacobi polynomials with parameters $(\alpha,\beta)$ 
defined by:
\begin{equation*}
\alpha =\frac{c}{2}(n-1)-1, \ \beta =\frac{c}{2}-1.
\end{equation*}
More precisely, the values of $(\alpha,\beta)$ are as follows:
\begin{equation*}
\begin{array}{|c|c|c|c|}
\hline 
 & c&\alpha&\beta \\
\hline
\R& 1&(n-3)/2 & -1/2 \\
\hline
\C & 2& n-2 & 0\\
\hline
\HH & 4 & 2n-3 & 1 \\
\hline
\OO & 8& 7& 3 \\
\hline
\end{array}
\end{equation*}
The  orthogonal measure associated to
the parameters $(\alpha,\beta)$ is $x^{\beta}(1-x)^{\alpha}$ over the
interval $[0,1]$. 
We generically
denote by $P_k(x)$ these polynomials, with $\deg(P_k)=k$ and
$P_k(1)=1$.
We let $\mu(x)$ denote their normalized associated orthogonal
measure
and $[P,Q]$ the corresponding scalar product on $\R[x]$ (so that 
$[P,Q]=\int P(x)Q(x)\mu(x)dx$ and $[1,1]=1$). Moreover, we have 
$[P_k,P_k]=d_k^{-1}$ where $d_k$ denotes the dimension of the
irreducible 
representation of $G$ associated to $P_k$ (e.g. when $X=\Sn$,
$d_k=\dim \Harm_k=\binom{n+k-1}{k}-\binom{n+k-3}{k-2}$).

The three-terms relation expresses $xP_k(x)$ as a linear
combination of the polynomials $P_i$:
\begin{equation*}
xP_k(x)=a_k P_{k+1}(x)+b_k P_{k}(x)+c_{k} P_{k-1}(x)
\end{equation*} 
for some sequences of rational numbers $(a_k)$, $(b_k)$, $(c_k)$. It
is enough for our purpose to know  that $(a_k)$ is bounded when $n$
and $k$ tend to $+\infty$ with $n/k$ tending to a finite limit. For
example,
when $X=\Sn$, 
\begin{equation*}
a_k=\frac{n-2+k}{n-2+2k}.
\end{equation*}

For all $(u,v)\in X$, we define
\begin{equation*}
t(u,v)=\left\{
\begin{array}{ll}(u\cdot v)=\cos\theta(u,v) &\text{if } X=\Sn \\
               \cos^2\theta(u,v) &\text{if } X=\Pn(K).
\end{array}
\right.
\end{equation*}
The zonal function on $X$ associated to $P_k$ is given by:
\begin{equation*}
(u,v)\to P_k(t(u,v)).
\end{equation*}

The so-called 'positivity property' related to these polynomials,
and underlying the linear programming method in $X$,
is the following: for all code $C\subset X$, and for
all $k\geq 0$,
\begin{equation*}
\sum_{u\in C}\sum_{v\in C} P_k(t(u, v))\geq 0.
\end{equation*}

\subsection{The linear programming method on $X^m$}
Now we consider the product spaces $X^m$. The positivity property
generalizes to the following:
\begin{lemma}\label{l1} Let $C\subset X^m$. Let us denote elements of
  $C$ by $u=(u_1,\dots,u_m)$ with $u_i\in X$.
For all $(k_1,\dots,k_m)\in \N^m$, 
\begin{equation*}
\sum_{u\in C}\sum_{v\in C} \prod_{i=1}^mP_{k_i}(t(u_i,v_i))\geq 0.
\end{equation*}
\end{lemma}
\begin{proof} This is the positivity property in the product space 
$X^m$. The group $G^m$ acts
transitively on $X^m$; the $G^m$-irreducible components
of $L^2(X^m)$ are the tensor products of the $G$-irreducible components
of each $L^2(X)$ and the associated zonal functions
are given by the polynomials in the $m$ variables $x_1,\dots,x_m$
\begin{equation*}
\prod_{i=1}^mP_{k_i}(x_i),\quad (k_1,\dots,k_m)\in \N^m
\end{equation*}
in the way:
\begin{equation*}
(u,v)\mapsto \prod_{i=1}^mP_{k_i}(t(u_i, v_i)).
\end{equation*}

\end{proof}

\begin{remark} In a sense, the polynomials $\prod_{i=1}^mP_{k_i}(x_i)$
are fake multivariate polynomials since the $m$ variables are
separated. As we shall see, for this reason it is much easier 
to deal with them, compared with the zonal polynomials for the
Grassmann space (see \cite{Bac06}).

\end{remark}

The polynomials $\prod_{i=1}^mP_{k_i}(x_i)$ generate the polynomial algebra\\ $\C[x_1,\dots,x_m]$, and are
orthogonal
for the product measure 
\begin{equation*}
\lambda \prod_{i=1}^m \mu(x_i)dx_i
\end{equation*}
with support $[-1,1]^m$ when $X=\Sn$, respectively $[0,1]^m$
otherwise, and  where $\lambda$ is chosen so that the total
measure is equal to $1$.  The associated scalar product on
$\R[x_1,\dots,x_m]$ is denoted by $[,]$. We take the following
notations: a multi-index in $\N^m$ is denoted by
$\uk=(k_1,\dots,k_m)$ and we define for $x=(x_1,\dots,x_m)$,
$P_{\uk}(x)=P_{\uk}(x_1,\dots,x_m):=\prod_{i=1}^mP_{k_i}(x_i)$, and 
$d_{\uk}:=\prod_{i=1}^m d_{k_i}$. 
Obviously we have,
for all $\uk$ and $\ul$,
\begin{equation}\label{eq scal prod}
\big[\Pk,\Pl\big]=\delta_{\uk,\ul} d_{\uk}^{-1}.
\end{equation}
Moreover, we define
$$\sigma(x):=\sum_{i=1}^m x_i.$$
For any angle $\theta$ we denote
\begin{equation*}
\left\{\begin{array}{ll}t=\cos\theta &\text{ if } X=\Sn\\
t=\cos^2\theta &\text{ if } X=\Pn(K)\\
\end{array}
\right.
\end{equation*}
Now we can formulate the usual associated linear programming bound:

\begin{proposition}\label{LP}
Assume $F\in \R[x_1,\dots,x_m]$ satisfies the conditions:
\begin{align*}
&(i) \ F=\sum_{\uk} f_{\uk} \Pk \text{ with }f_{\uk}\geq 0 \text{  for all }\uk,
\text{  and }f_0>0\\
&(ii)
\left\{\begin{array}{lll} \text{If } X=\Sn, &F(x_1,\dots,x_m)\leq 0
&\text{for all }(x_1,\dots,x_m)\in [-1,1]^m \\
&&\text{such that
}\sigma(x)\leq m\cos\te=mt\\ 
\text{If } X\neq\Sn, &F(x_1,\dots,x_m)\leq 0
&\text{for all }(x_1,\dots,x_m)\in [0,1]^m \\
&&\text{such that }\sigma(x)\leq m\cos^2\te=mt
\end{array}\right.
\end{align*}

Then, any code $C$ in $X^m$ with minimum angular distance $\te$ satisfies 
\begin{equation*}
|C| \leq \frac{F(1,\dots,1)}{f_0}
\end{equation*}
\end{proposition}

\begin{proof} We reproduce the standard argument. Let 
\begin{equation*}
S:=\sum_{u\in C}\sum_{v\in C} F(t(u_1, v_1),\dots,t(u_m, v_m)).
\end{equation*}
The pairs $(u,v)$ with $u=v$ contribute in this sum for
$|C|F(1,\dots,1)$. From condition (ii) and the assumption
that for $u\neq v\in C$, $\cos\theta(u,v)\leq\cos\te$, the other terms
are non positive. Hence, $S\leq |C|F(1,\dots,1)$.

On the other hand, we have 
$$S=\sum_{\uk} f_{\uk}\big(\sum_{u, v\in C}
\Pk(t(u_1,v_1),\dots,t(u_m,v_m)).$$
The term corresponding to $\uk=0=(0,\dots,0)$ gives $f_0|C|^2$ while
the other terms are non-negative from the positivity property of the
polynomials $\Pk$ (Lemma \ref{l1}). Hence $S\geq f_0|C|^2$.
The two inequalities lead to the announced bound.

\end{proof}

\subsection{Examples of small degree }
Let us work out the case of polynomials of small degree. 

\begin{enumerate} 
\item $X=\Sn$
\begin{enumerate}
\item Degree $1$: we take $F=(x_1+\dots+x_m)-mt$. Since
$P_1(x)=x$,
$F$ satisfies the hypo\-the\-sis of Proposition \ref{LP} if and only if $t<0$.
We obtain:
\begin{equation}\label{deg1}
\text{If } \cos\theta=t<0,\quad |C|\leq 1-\frac{1}{t}.
\end{equation}
\item Degree $2$: we take
$F=\big((x_1+\dots+x_m)+m\big)\big((x_1+\dots+x_m)-mt).$
We have 
\begin{align*}
F&=(x_1+\dots+x_m)^2+m(1-t)(x_1+\dots+x_m)-m^2t\\
&=\sum x_i^2+2\sum_{i<j} x_ix_j +m(1-t)\sum x_i -m^2t\\
&=\sum (x_i^2-\frac{1}{n})+2\sum_{i<j} x_ix_j +m(1-t)\sum x_i +\frac{m}{n}-m^2t\\
\end{align*}

Since $P_2(x)=(x^2-1/n)/(1-1/n)$, $F$ satisfies the hy\-po\-the\-sis of
Proposition \ref{LP} if and only if $\frac{m}{n}-m^2t>0$.
We obtain:
\begin{equation}\label{deg2}
\text{If } \cos\theta=t<\frac{1}{mn},\quad |C|\leq \frac{2mn(1-t)}{1-mnt}.
\end{equation}
The two bounds take the value $1+mn$
at their cros\-sing point, cor\-res\-ponding  to $t=-1/mn$.
.
\end{enumerate}

\item $X=\Pn(K)$, Degree $1$: we have, up to a multiplicative factor, $P_1(x)=x-\frac{\beta+1}{\alpha+\beta+2}=x-\frac{1}{n}$.
We take
$F=(x_1+\dots+x_m)-mt=(x_1-\frac{1}{n})+\dots+(x_m-\frac{1}{n})+\frac{m}{n}-mt$.
$F$ satisfies the hypothesis of Proposition \ref{LP} if and only if
$t<1/n$. We obtain:
\begin{equation}\label{deg1 Pn}
\text{If } \cos^2\theta=t<\frac{1}{n},\quad |C|\leq \frac{1-t}{1/n-t}.
\end{equation}

\end{enumerate}

\subsection{Christoffel-Darboux formula and an explicit bound}
It remains to apply
the standard method with Christoffel-Darboux formula. For $\uk=(k_1,\ldots,k_m)$ and $\ul=(l_1,\ldots,l_m)$,
the notation $\ul\leq \uk$ stands for: $l_i\leq k_i$ for all 
$1\leq i\leq m$.

\begin{proposition} \label{prop CDF}
Let $y=(y_1,\dots,y_m)\in \R^m$ and $\uk=(k_1,\ldots,k_m)\in\N^m$, and define
\begin{equation*}
K_{\uk}(x,y):=\sum_{\ul\leq \uk}d_{\ul}\Pl(x)\Pl(y)
=\prod_{j=1}^m \big(\sum_{i=0}^{k_j}d_iP_i(x_j)P_i(y_j)\big)
\end{equation*}
and
\begin{equation*}
N_{\uk}(x,y):=\sum_{t=1}^m d_{k_t}a_{k_t}
Q_{k_t}(x_t,y_t)
\prod_{j\neq t} \big(\sum_{i=0}^{k_j}d_iP_i(x_j)P_i(y_j)\big)
\end{equation*}
where
\begin{equation*}
Q_{k_t}(x_t,y_t):=P_{k_t+1}(x_t)P_{k_{t}}(y_t)-P_{k_{t}}(x_t)P_{k_t+1}(y_t).
\end{equation*}
Then we have the Christoffel-Darboux type formula:
\begin{equation*}
K_{\uk}(x,y)=\frac{N_{\uk}(x,y)}{\sigma(x)-\sigma(y)}.
\end{equation*}

\end{proposition}

\medskip

\begin{proof} Since $\sigma(x)-\sigma(y)=\sum_{t=1}^m x_t-\sum_{t=1}^m y_t
=\sum_{t=1}^m (x_t-y_t)$,
\begin{align*}
(\sigma(x)-\sigma(y))&K_{\uk}(x,y)=\Big(\sum_{t=1}^m (x_t-y_t)\Big)
\prod_{j=1}^m \big(\sum_{i=0}^{k_j}d_iP_i(x_j)P_i(y_j)\big)\\
&=\sum_{t=1}^m \Big( (x_t-y_t)\sum_{i=0}^{k_t}d_iP_i(x_t)P_i(y_t)\Big)
\prod_{j\neq t} \big(\sum_{i=0}^{k_j}d_iP_i(x_j)P_i(y_j)\big)
\end{align*}
The Christoffel-Darboux formula for the polynomials $P_k$ gives:
\begin{equation*}
(x_t-y_t)\sum_{i=0}^{k_t}d_iP_i(x_t)P_i(y_t)=d_{k_t}a_{k_t}Q_{k_t}(x_t,y_t)
\end{equation*}
with the notations of the proposition, hence the result.

\end{proof}

Following the standard method, we apply Proposition \ref{LP}
to the function 
$$\frac{N_{\uk}(x,y)^2}{\sigma(x)-\sigma(y)}.$$

\begin{proposition} \label{prop nabound}
Let $y=(y_1,\dots,y_m)\in \R^m$ and $\uk\in\N^m$, and define

\begin{equation*}
F(x):=\frac{N_{\uk}(x,y)^2}{\sigma(x)-\sigma(y)}=(\sigma(x)-\sigma(y))K_{\uk}(x,y)^2=K_{\uk}(x,y)N_{\uk}(x,y).
\end{equation*}

If $y$ satisfies the conditions:
\begin{enumerate}
\item $P_i(y_t)\geq 0$ for all $0\leq i\leq k_t$ and for all $1\leq t\leq m$
\item $P_{k_t+1}(y_t)\leq 0$ for all $1\leq t\leq m$
\end{enumerate}
then $F$ satisfies the hypothesis of Proposition \ref{LP} for all
$\te$ such that $mt\leq\sigma(y)$.
Consequently, for any code $C$ in $X^m$ with minimum angular distance $\te$,
\begin{equation}\label{bb}
|C|\leq \frac{\Big(\sum_{t=1}^m
d_{k_t}a_{k_t}\big(P_{k_t}(y_t)-P_{k_t+1}(y_t)\big)\prod_{j\neq t}\big(\sum_{i=0}^{k_j}d_iP_i(y_j)\big)\Big)^2}
{-(m-\sigma(y))\sum_{t=1}^m d_{k_t} a_{k_t} P_{k_t}(y_t)P_{k_t+1}(y_t) \prod_{j\neq t}\left(\sum_{i=0}^{k_j} d_i(P_i(y_j))^2\right)}.
\end{equation}
\end{proposition}

\medskip

\begin{proof} Clearly, under the assumptions (i) and (ii), $K_{\uk}(x,y)$ and $N_{\uk}(x,y)$ have
non-negative coefficients on the $\Pl$. This is enough to ensure
that it is also the case for the product $K_{\uk}(x,y)N_{\uk}(x,y)$
(recall that the product of two polynomials with non-negative coefficients on 
the $P_k$ also has non-negative coefficients on the $P_k$. This
property transfers straightforwardly to the $\Pk$; it is anyway
general
to any family of zonal polynomials).

Obviously the sign of $F(x)$ is the sign of
$\sigma(x)-\sigma(y)$ so the conditions of
Proposition \ref{LP} are fulfilled.

It remains to compute $f_0=[F,1]$ and $F(1,\dots,1)$. 
\begin{align*}
&[F,1]=[K,N]\\
& = \Big[\prod_{j=1}^m\big(\sum_{i=0}^{k_j}d_i P_i(x_j)P_i(y_j)\big),
\sum_{t=1}^m d_{k_t} a_{k_t} Q_{k_t}(x_t,y_t) \prod_{j\neq t} \big(\sum_{i=0}^{k_j}d_i P_i(x_j)P_i(y_j)\big)\Big] \\
& = \sum_{t=1}^m d_{k_t} a_{k_t} \Big[\prod_{j=1}^m\big(\sum_{i=0}^{k_j}d_i P_i(x_j)P_i(y_j)\big),
 Q_{k_t}(x_t,y_t) \prod_{j\neq t} \big(\sum_{i=0}^{k_j}d_i P_i(x_j)P_i(y_j)\big)\Big]  \\
& = \sum_{t=1}^m d_{k_t} a_{k_t} \Big[\sum_{i=0}^{k_t}d_i
    P_i(x_t)P_i(y_t),Q_{k_t}(x_t,y_t)\Big]\\
&\qquad\qquad\qquad\qquad\qquad\qquad.\prod_{j\neq t}\Big[\sum_{i=0}^{k_j}d_i P_i(x_j)P_i(y_j),\sum_{i=0}^{k_j}d_i P_i(x_j)P_i(y_j)\Big]  \\
& = \sum_{t=1}^m d_{k_t} a_{k_t} \big(-P_{k_t}(y_t)P_{k_t+1}(y_t)\big)\prod_{j\neq t}\big(\sum_{i=0}^{k_j}d_i (P_i(y_j))^2\big)
\end{align*}
where the last equality follows from \eqref{eq scal prod}.

Let us now compute $F(1,\dots,1)$. We have:
\begin{equation*}
F(1,\dots,1)=\frac{N_{\uk}(\underline{1},y)^2}{m-\sigma(y)}
\end{equation*}
and
\begin{equation*}
N_{\uk}(\underline{1},y)=\sum_{t=1}^m
d_{k_t}a_{k_t}\big(P_{k_t}(y_t)-P_{k_t+1}(y_t)\big)\prod_{j\neq t}\big(\sum_{i=0}^{k_j}d_iP_i(y_j)\big).
\end{equation*}
Applying the resulting bound of Proposition \ref{LP} leads to the
announced bound.
\end{proof}

We proceed to choose the parameters $y$ and $\uk$ such that the
conditions of Proposition \ref{prop nabound} will be satisfied. We
follow the standard method. We first choose the multi-index $\uk$ such
that $mt\leq\sum_{t=1}^m z_{k_t}$, where $z_{k_t}$ is the largest zero
of $P_{k_t}$. The interlacing property of the zeros of the orthogonal
polynomials $\Pk$ guarantees that there exists $y$ such that $z_{k_t}\leq y_{t}\leq z_{k_t+1}$ and
\begin{equation*}
P_{k_t}(y_t)+P_{k_t+1}(y_t)=0
\end{equation*}
for all $1\leq t\leq m$. Thus, $P_i(y_t)>0$ and $P_{k_t+1}(y_t)<0$ for all $0\leq i\leq k_t$ and $1\leq t\leq m$. 

Now we have:
\begin{align*}
f_0=[F,1]&= \sum_{t=1}^m d_{k_t}a_{k_t}\big(P_{k_t}(y_t)\big)^2\prod_{j\neq
  t}\big(\sum_{i=0}^{k_j}d_i (P_i(y_i))^2\big)\\
&=\sum_{t=1}^m a_{k_t}\sum_{\substack{\ul\leq \uk\\l_t=k_t}}
  d_{\ul}\big(P_{\ul}(y)\big)^2
\end{align*}
and
\begin{align*}
N_{\uk}(\underline{1},y)&=2\sum_{t=1}^m
d_{k_t}a_{k_t}P_{k_t}(y_t)\prod_{j\neq
  t}\big(\sum_{i=0}^{k_j}d_iP_i(y_j)\big)\\
&=2\sum_{t=1}^m a_{k_t}\sum_{\substack{\ul\leq \uk\\l_t=k_t}}
  d_{\ul}P_{\ul}(y).
\end{align*}
With Cauchy-Schwartz inequality (applied twice),
\begin{align*}
F(1,\dots,1)&= \frac{4}{(m-\sigma(y))}\Big(\sum_{t=1}^m a_{k_t}\sum_{\substack{\ul\leq \uk\\l_t=k_t}}
  d_{\ul}P_{\ul}(y)\Big)^2\\
&\leq \frac{4}{(m-\sigma(y))}\Big(\sum_{t=1}^m a_{k_t}\Big)
\Big(\sum_{t=1}^m a_{k_t}\big(\sum_{\substack{\ul\leq \uk\\l_t=k_t}}
  d_{\ul}P_{\ul}(y)\big)^2\Big)\\
&\leq \frac{4}{(m-\sigma(y))}\Big(\sum_{t=1}^m a_{k_t}\Big)
\Big(\sum_{t=1}^m a_{k_t}\big(\sum_{\substack{\ul\leq \uk\\l_t=k_t}}d_{\ul}\big)
\big(\sum_{\substack{\ul\leq \uk\\l_t=k_t}}
  d_{\ul}\big(P_{\ul}(y)\big)^2\big)\Big)\\
&\leq \frac{4}{(m-\sigma(y))}\Big(\sum_{t=1}^m a_{k_t}\Big)
\Big(\sum_{\ul\leq \uk}d_{\ul}\Big)\Big(\sum_{t=1}^m a_{k_t}
\sum_{\substack{\ul\leq \uk\\l_t=k_t}}
  d_{\ul}\big(P_{\ul}(y)\big)^2\Big)\\
&= \frac{4}{(m-\sigma(y))}\Big(\sum_{t=1}^m a_{k_t}\Big)
\prod_{t=1}^m \big(\sum_{i=0}^{k_t}d_i\big)f_0
\end{align*}
Denote $D_{k_t}:=\sum_{i=0}^{k_t}d_i$. We obtain 
\begin{equation*}
|C| \leq \frac{4\big(\sum_{t=1}^m a_{k_t}\big)\prod_{t=1}^m
 D_{k_t}}{m-\sigma(y)}.
\end{equation*}

We summarize the above result in the following statement:

\begin{proposition}
For any code $C$ in $X^m$ with minimum angular distance $\theta$, for
any multi-index $\uk$ such that $mt\leq \sum_{t=1}^m z_{k_t}$, let
$y_t$ satisfy $z_{k_t}\leq y_{t}\leq z_{k_t+1}$ and
\begin{equation*}
P_{k_t}(y_t)+P_{k_t+1}(y_t)=0,
\end{equation*}
then 
\begin{equation}\label{bound LP}
|C| \leq \frac{4\big(\sum_{t=1}^m a_{k_t}\big)\prod_{t=1}^m
 D_{k_t}}{m-\sigma(y)}.
\end{equation}

\end{proposition}

\begin{remark}
Using the so-called adjacent polynomials instead of the Gegenbauer polynomials, an enhancement of \eqref{bound LP} was derived for $m=1$ \cite{Lev98, EZ01}. It seems that this can be generalized for all $m$.

\end{remark}

\subsection{A bound for the asymptotic rate}

Now we consider the limit when $n\to+\infty$ of the rate $R(C):=\frac{1}{n}\ln|C|$ (of
course the space $\PO$ is not concerned anymore) of the codes $C$ of
$X^m$. We derive an  upper bound for this limit from (\ref{bound
  LP}). The next theorem settles the result obtained that way only in
the case $X=\Pn(K)$ because this bound,
in the case of $X=\Sn$, turns out to be the same as the one obtained from
the trivial isometric embedding $\Snm\to S^{mn-1}$ (see Remark \ref{big remark}).

\begin{theorem} \label{bound} 
Let $C$ be a code in $X^m$, $X=\Pn(K)$, with minimum angular distance
$\theta$, and let $(\theta_1,\ldots,\theta_m)\in [0,\pi/2]^m$ satisfy $\sum_{t=1}^m \cos^2 \theta_t=m\cos^2\theta$. Then, when $n\to\infty$,
\begin{equation}    \label{borne}
R(C)\lesssim c(R_{LP}(\te_1)+\ldots+R_{LP}(\te_m)),
\end{equation}
where $R_{LP}$ is defined in \eqref{theo lp sphere}.
\end{theorem}

\begin{proof} Same as in \cite{KL78}, 
involving the asymptotic estimate of $z_k$. We reproduce it here: Consider an infinite sequence $k(n)$ such that $2k(n)/cn$ tends to a finite limit $\rho$ as $n$ tends to infinity. Then \cite{KL78}
\begin{equation*}
\lim_{n\to\infty}z_{k(n)}= 4\frac{\rho^{-1}+1}{(\rho^{-1}+2)^2}
\end{equation*}
and, since from \cite{Lev98},
\begin{equation*}
D_k\simeq \binom{\frac{c}{2}n+k-1}{k}^2
\end{equation*}
we have
\begin{equation*}
\lim_{n\to\infty}\frac{1}{n}\ln D_{k(n)}= \lim_{n\to\infty}\frac{2}{n}\ln\binom{\frac{c}{2}n+k(n)-1}{k(n)}=c\big((1+\rho)\ln(1+\rho)-\rho\ln\rho\big).
\end{equation*}

Inverting the conditions
\begin{equation*}
\cos^2\theta_t=4\frac{\rho_t^{-1}+1}{(\rho_t^{-1}+2)^2}
\end{equation*}
leads to
\begin{equation*}
\rho_t=\frac{1-\sin\te_t}{2\sin\te_t}
\end{equation*}

Let $k_t=\lfloor\rho_t n\rfloor$, and let $y_t$ satisfy $z_{k_t}\leq
y_t\leq z_{k_t+1}$ and 
$P_{k_t}(y_t)+P_{k_t+1}(y_t)=0$ (the existence of $y_t$ is guaranteed
by the interlacing property of the zeros of the Jacobi
polynomials). Then  from \eqref{bound LP},
\begin{equation*}
|C| \leq \frac{4\big(\sum_{t=1}^m a_{k_t}\big)\prod_{t=1}^m
 D_{k_t}}{m-\sigma(y)}.
\end{equation*}

Since $\sigma(y)\simeq m\cos^2\theta$ and the expression
$\frac{4\sum_{t=1}^m a_{k_t}}{m-\sigma(y)}$ 
has a finite limit when $k_t/n$
tends to $\rho_t$,   the rate $R(C)$ satisfies
\begin{equation*}
R(C)\lesssim \frac{1}{n}\sum_{t=1}^m\ln D_{k_t}\simeq
\sum_{t=1}^m c\big((1+\rho_t)\ln(1+\rho_t)-\rho_t\ln(\rho_t)\big).
\end{equation*}
\end{proof}

\begin{remark}  \label{big remark}
\begin{itemize}

\item It is worth noticing that the choice $\theta_t=\theta$ in
  \eqref{borne} yields to the bound
\begin{equation}\label{borne2}
R(C) \lesssim cm R_{LP}(\te).
\end{equation}
This bound can be derived more easily, since every code in $X^m$ is
also a code in the $cmn$-th dimensional unit sphere
(combining the mapping $\beta$ for $m=1$ and the  obvious mapping $\Snm\to S^{mn-1}$).
It turns out that, since the function $R_{LP}(\theta)$ {\em as a function of $t=\cos^2\theta$} is not convex, the bound \eqref{borne}
slightly improves on \eqref{borne2}. We discuss this in more details in
the next subsection.

\item The same method applied to $X=\Sn$ would lead to:
\begin{equation*}
R(C)\lesssim R_{LP}(\te_1)+\ldots+R_{LP}(\te_m), \text{ for all }
\theta_t \text{ such that } \sum_{t=1}^m \cos \theta_t=m\cos\theta.
\end{equation*} 
But the function $R_{LP}(\theta)$ {\em as a function of $t=\cos\theta$}
is convex, therefore the choice of
$(\theta_1,\dots,\theta_m)$ that minimizes the right hand side is
$\theta_1=\dots=\theta_m=\theta$, yielding \eqref{borne2}.
\end{itemize}
\end{remark}

\subsection{Analysis of \eqref{borne} versus \eqref{borne2}}

Let $C^2$ be the set of continuous, twice differentiable functions with continuous second derivative. For a function $f$ defined on
$[0,1[$, of class $C^2$, we denote:
\begin{equation*}
f^{(m)}(t):= \min_{\substack{t_1,\dots,t_m\in [0,1[\\\sum_{i=1}^m
    t_i=mt}}\frac{f(t_1)+\dots+f(t_m)}{m}.
\end{equation*}
Clearly, if $f$ is convex on $[0,1[$, we have $f^{(m)}=f$, and, if
$f\leq g$, $f^{(m)}\leq g^{(m)}$. It is also easy to see that
$f^{(m')}\leq f^{(m)}$ when $m$ divides $m'$.

The function we are interested in is $f(t)=R_{LP}(\theta)$ where
$t=\cos^2\theta$. We have
\begin{equation*}
f(t)=(1+\rho(t))\ln(1+\rho(t))-\rho(t)\ln(\rho(t))
\end{equation*}
where
\begin{equation*}
\rho(t)=\frac{1}{2}\big(-1+(1-t)^{-1/2}\big).
\end{equation*}
One can check that the second derivative of $f$ takes negative values
on some interval $[0,t_0]$, $t_0\simeq 0.208$, and then takes positive values on $[t_0,1[$.
The function $f$ is an increasing function, with $f(0)=0$,  first concave then convex.
We consider the function $g$ on $[0,1[$, whose graph ${\mathcal C}_g$ determines the
convex hull of the portion of plane above the graph ${\mathcal C}_f$ of $f$. The
function $g$ is uniquely determined by the conditions:
\begin{equation*}
\left\{
\begin{array}{l}
g\leq f\\
g\text{ is convex }\\
g \text{ is maximal with these properties }
\end{array}
\right.
\end{equation*}
Let us denote by $t_1$ the unique value for which the tangent at
$(t_1,f(t_1))$  to  ${\mathcal C}_f$ contains the origin $(0,0)$. The value
$t_1\simeq 0.379$
is the unique solution to 
$$f(t)=f'(t)t$$
and the slope of the tangent to ${\mathcal C}_f$ at $t_1$ equals $f'(t_1)\simeq 1.089$.
Then the function $g$ is defined by:
\begin{equation*}
\left\{
\begin{array}{ll}
g(t)=f'(t_1)t\simeq 1.089t &\text{ for all }t\in [0,t_1]\\
g(t)=f(t) &\text{ for all }t\in [t_1,1[\\
\end{array}
\right.
\end{equation*}
Since $g$ is convex, we have for all $m$ and all $t\in [0,1[$,
$g^{(m)}(t)=g(t)\leq f^{(m)}(t)$. In other words, on $[0,t_1]$,
$f^{(m)}$ is somewhere between $g$ and $f$, and on $[t_1,1[$,
$f^{(m)}=f$.
Clearly, when $m\to+\infty$, $f^{(m)}\to g$. Also, the maximum
$\delta$ of
$f(t)-g(t)$
is an upper bound for the maximum of $f(t)-f^{(m)}(t)$. Numerical
calculation gives $\delta\simeq 0.016$. Considering our primary goal, i.e., to compare \eqref{borne} and \eqref{borne2}, this means that the improvement of \eqref{borne2} upon \eqref{borne} is upper-bounded by $0.016m$.

It seems difficult to determine the optimal choice of
$(t_1,\dots,t_m)$ that minimizes the quotient 
$\frac{f(t_1)+\dots+f(t_m)}{m}$. A natural choice is
$(t_1,\dots,t_m)=(0,0,..,mt/r,\dots,mt/r)$ with $r$ non-zero and equal
coordinates. In that case,
$\frac{f(t_1)+\dots+f(t_m)}{m}=\frac{r}{m}f(\frac{mt}{r})$
and requires $t<r/m$. If $t=\frac{rt_1}{m}$, it is certainly the best
choice since then the resulting point lies on ${\mathcal C}_g$. 
Numerical experiments seem to show that, for $m=2,3$, and $t<1/m$,
$r=1$ does minimize the quotient 
$\frac{f(t_1)+\dots+f(t_m)}{m}$.

\section{Bounds for codes in the Grassmann and Stiefel
  manifolds}\label{sec grass stief}

In this section, we summarize the consequences of the above results
for Grassmann and Stiefel codes. Following a standard notation in
coding theory, we denote by $A(X,d)$, the maximal number of elements
of a code $C$ of the space $X$ with minimum distance $d$.

We have proved in the section \ref{subsec C} that the size of Grassmannian codes with minimal chordal
distance $d=\sqrt{m-s}$ is upper bounded by the size of codes in
$\Pn(K)^m$ with minimal angular distance $\theta$, where
$\cos^2\theta=s/m$. Thus we have proved that:
\begin{equation}\label{A bound grass}
A(\Gmn(K),d)\leq A(\Pn(K)^m,\theta)\text{ with } \theta=\arccos{\sqrt{1-d^2/m}}.
\end{equation}
Linear programming bounds on $A(\Pn(K)^m,\theta)$ were derived in Section \ref{sec products}.

We believe that these bounds are not good in general for finite values
  of the parameters, because we use only a rough estimate of
  $\sigma(p,q)=\trace(\pi_p\circ\pi_q)$ in the inequality
  \eqref{ineq2} (we replace $\sum_{1\leq i,j\leq m} (e_i\cdot e'_j)^2$
  with $\sum_{1\leq i\leq m} (e_i\cdot e'_i)^2$). If we compare the bounds obtained with the zonal
  polynomials of small degree, \eqref{deg1 Pn}
  is worse than the simplex bound, obtained from the zonal
  polynomial of degree $1$ of
  $\Gmn(\R)$.
Moreover, numerical experiments for small parameters $m$ and $n$
  (with the package LRS, by David Avis,
  http://cgm.cs.mcgill.ca/\verb+~+avis/C/lrs.html),  confirms that the bounds
  obtained from the zonal polynomials of $\Gmn(\R)$ are sharper than
  the ones obtained from Proposition \ref{LP} for $X=\Pn(\R)^m$.

\begin{figure}
\begin{center}
\includegraphics*{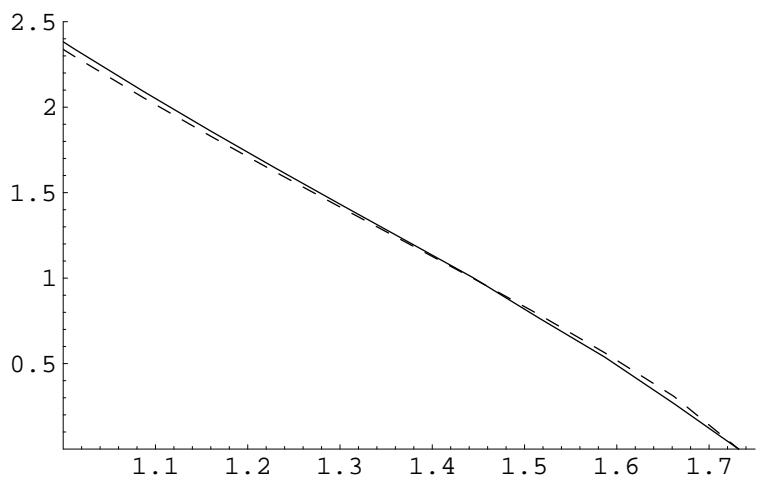}
\caption{Upper bounds on the asymptotic rate of real Grassmannian codes with minimum chordal distance $d$, $m=3$. The solid line is $R_1(d)$, and the dashed line is $R_2(d)$ (see Theorem \ref{bound grass})}    \label{fig comp}
\end{center}
\end{figure}

Surprisingly, the consideration of $\Pn(K)^m$ allows us to obtain
better bounds for the asymptotic rate than the ones obtained
previously by either the isometric embedding given in \cite{CHS} of
$\Gmn$ into a unit sphere of the dimension $(n-1)(n+2)/2$ (see also
Remark \ref{remark}), or the spectral method
developed in \cite{Bac06} with the zonal polynomials of $\Gmn$. 
We summarize the new bound we have obtained in the next theorem:

\begin{theorem} \label{bound grass final} 
Let $C$ be a code in $\Gmn(K)$ with minimal chordal distance $d=\sqrt{m-s}$,
and let $\theta=\arccos(\sqrt{s/m})$. Then, when $n\to +\infty$,
\begin{equation} \label{bound2 gmn}
R(C)\lesssim\min\left\{R_1(d),R_2(d)\right\},
\end{equation}
where 
\begin{equation}
R_1(d) = \min_{\substack{(\theta_1,\dots,\theta_m)\in
    [0,\pi/2]^m\\ \sum_{i=1}^m\cos^2\theta_i=m\cos^2\theta}}c\big(R_{LP}(\theta_1)+\dots +R_{LP}(\theta_m)\big)
\end{equation}
and
\begin{equation}
R_2(d) = cmR_S(\te).
\end{equation}
\end{theorem}
The bounds $R_1(d)$ and $R_2(d)$ are depicted in Figure \ref{fig comp}.

\smallskip
We have proved in Section \ref{subsec C} that the size of Stiefel codes with minimal chordal
distance $d=\sqrt{2}\sqrt{m-s}$ is upper-bounded by the size of codes in
$\big(S^{cn-1}\big)^m$ with minimal angular distance $\theta$, where
$\cos\theta=s/m$. In Section \ref{sec products}, we derived linear programming bounds for codes in $\big(S^{cn-1}\big)^m$, thus implying bounds for Stiefel codes.
 These bounds are, up to our knowledge, the first general
bounds for Stiefel codes, and we believe that they are rather sharp.
For the asymptotic rate, the best result is obtained in \eqref{bound1 stiefel}.
We summarize these results in the following theorem:

\begin{theorem} \label{bound stiefel final} 
With the previous notations, and $\theta=\arccos(1-d^2/2m)$,
\begin{enumerate}
\item 
$A(\Vmn(K),d)\leq A(\big(S^{cn-1}\big)^m,\theta)$
\item When $n\to+\infty$, $R(C)\lesssim cmR_S(\theta)$
\end{enumerate}
\end{theorem}

\section{Conclusions}   \label{sec con}

Using relations between Grassmann and Stiefel manifolds and other spa\-ces, we derive new bounds on the size of Grassmannian codes (Theorem \ref{bound grass final}) and Stiefel codes (Theorem \ref{bound stiefel final}). These are the best known asymptotic bounds on the rate of Grassmannian and Stiefel codes.
\smallskip

\newpage
\begin{center}
{\bf Acknowledgement}
\end{center}

We thank Jean Creignou for helpful discussions about spherical and Grassmannian codes.

\end{document}